\documentclass[12pt]{article}
\usepackage{amsfonts}
\usepackage{amssymb}
\usepackage{amsmath}
\usepackage{amsthm}
\usepackage{verbatim}

\newtheorem{thm}{Theorem}
\newtheorem{lem}{Lemma}

\theoremstyle{definition} 
\newtheorem{defn}{Definition}

\newcommand{\p}{\varphi}

\newcommand{\pa}{\varphi^t}

\newcommand{\Om}{\Omega}
\newcommand{\om}{\omega}

\newcommand{\R}{\mathbb{R}}

\newcommand{\N}{\mathbb{N}}
\newcommand{\T}{\mathbb{T}}

\newcommand{\D}{\mathcal{D}}

\newcommand{\bmo}{{\rm BMO}}

\newcommand{\blo}{{\rm BLO}}

\newcommand{\sumA}{\sum_{n:2^{-n}<|Q|}}
\newcommand{\sumB}{\sum_{n:2^{-n}\ge|Q|}}

\newcommand{\intav}{-\!\!\!\!\!\!\int}
\newcommand{\intq}{\frac{1}{|Q|} \, \int_Q \,}

\DeclareMathOperator*{\essinf}{ess\ inf}
\DeclareMathOperator*{\esssup}{ess\ sup}

\def\XXint#1#2#3{{\setbox0=\hbox{$#1{#2#3}{\int}$}
     \vcenter{\hbox{$#2#3$}}\kern-.5\wd0}}

\begin{document}

\title{Geometric-arithmetic averaging of dyadic
    weights}

\renewcommand{\thefootnote}{\fnsymbol{footnote}}

\author{Jill Pipher\footnote{Research supported by the NSF under grant
    number DMS0901139.}\\
Department of Mathematics\\
Brown University\\
Providence, RI 02912\\
{\tt jpipher@math.brown.edu}
\and Lesley A.~Ward\\ 
School of Mathematics and Statistics\\
University of South Australia\\
Mawson Lakes, SA 5095\\
Australia\\
{\tt lesley.ward@unisa.edu.au} \and Xiao Xiao$^\dag$ \\
Department of Mathematics\\
Brown University\\
Providence, RI 02912\\
{\tt xxiao@math.brown.edu}}

\renewcommand{\thefootnote}{\arabic{footnote}}


\maketitle


\begin{abstract}
The theory of (Muckenhoupt) weights arises in many areas of
analysis, for example in connection with bounds for singular
integrals and maximal functions on weighted spaces. We prove
that a certain averaging process gives a method for
constructing $A_p$~weights from a measurably varying family of
dyadic $A_p$ weights. This averaging process is suggested by
the relationship between the $A_p$~weight class and the space
of functions of bounded mean oscillation. The same averaging
process also constructs weights satisfying reverse H\"older
($RH_p$) conditions from families of dyadic $RH_p$ weights, and
extends to the polydisc as well.
\end{abstract}


\section{Introduction}
\label{sec:introduction} Several classes of functions are
defined in terms of a property that the function must satisfy
on each interval, with a uniform constant. Well
known examples from harmonic analysis and complex analysis
include Muckenhoupt's $A_p$ weights, the reverse-H\"older
weight classes~$RH_p$, the class of doubling weights, and the
space \bmo\ of functions of bounded mean oscillation. Such
classes have strictly larger dyadic analogues, where the
defining property is required to hold only on the dyadic
intervals. Certain types of averaging provide a
bridge between these dyadic counterparts and the original
function classes. Specifically, these averages convert each
suitable family of functions in the dyadic class to a single
function in the smaller, nondyadic class. We can think of
averaging as an improving operation, in this sense.

An easily stated example is the following. If all translates of
a function~$f$ defined on the unit circle $\T = [0,1]$ are in
dyadic~\bmo, or equivalently if $f$ is in dyadic~\bmo\ on every
translated grid of dyadic intervals on the circle, then the
function~$f$ itself is in true~\bmo. This result is a special
case of a theorem in~\cite{GJ}, applied to the identity
\[
  f(x)
  = \int_0^1 \tau_t f(x + t) \, dt,
\]
where for $t\in\R$ the translation operator $\tau_t$ is defined
by $\tau_t f(\cdot) := f(\cdot - t)$, and $x + t$ is to be
interpreted as $x + t \mod 1$.

Now, what if a function $f$ on~$\T$ can be written as the
\emph{translation-average}
\[
  f(x)
  := \int_0^1 f^t(x + t) \, dt
\]
of dyadic \bmo\ functions $\{f^t\}_{t\in[0,1]}$ that are
\emph{not} identical translates of each other? If they satisfy
the hypotheses of~\cite{GJ}, then still $f$ is in true~\bmo.
However, the analogous statements can fail for $A_p$~weights,
for $RH_p$~weights, and for doubling weights~\cite{W}.

In this paper we show that a different type of averaging works
for both $A_p$ and~$RH_p$ (Theorems~\ref{thm:Apexplog}
and~\ref{thm:RHpexplog}). This is the
\emph{geometric-arithmetic average} defined by
\[
  \Om(x)
  := \exp \left\{\int_0^1 \log \om^t(x + t) \, dt\right\},
\]
where $\{\om^t\}_{t\in [0,1]}$ is a suitable family of weights
in $A_p^d$ or~$RH_p^d$.

We also observe that translation-averaging does work for both
$A_p$ and~$RH_p$ under the additional assumption that the
functions~$\om^t$ are doubling weights, not just dyadic
doubling weights (Theorem~\ref{thm:doubling}).

All these results generalize to the polydisc
(Theorems~\ref{thm:productApRHpexplog}
and~\ref{thm:Bidoubling}).

The paper is organized as follows. In Section~\ref{sec:Thms},
we state our geometric-arithmetic averaging results on the
circle. In Section~\ref{sec:definitions}, we collect the
definitions and background results used in the paper. Also,
Lemma~\ref{lem:ApRHpcharn} in that section gives a unified
characterization of weights in~$A_p$ for $1 \leq p \leq
\infty$, $RH_p$ for $1 < p \leq \infty$, and their dyadic
counterparts, in terms of conditions on the oscillation of
their logarithms. We take some care throughout in tracing the
dependence of the various constants. In
Section~\ref{sec:explogApRHp}, we prove geometric-arithmetic
averaging for $A_p^d$ and $RH_p^d$ weights
(Theorems~\ref{thm:Apexplog} and~\ref{thm:RHpexplog}). In
Section~\ref{sec:dbltranslproof}, we prove
translation-averaging for $A_p^d$ and $RH_p^d$ weights that are
doubling (Theorem~\ref{thm:doubling}). In
Section~\ref{sec:polydisc}, we generalize our results to the
polydisc (Theorems~\ref{thm:productApRHpexplog}
and~\ref{thm:Bidoubling}).

\section{Geometric-arithmetic averaging on the circle}
\label{sec:Thms}

In~\cite{W}, examples are constructed to show that, for
$\{\om^t\}_{t\in[0,1]}$ a measurably varying family of dyadic
$A_p^d$ weights (for arbitrary $p$ with $1 \leq p \leq \infty$)
or dyadic $RH_p^d$ weights (for arbitrary $p$ with $1 < p <
\infty$) on the circle $\T = [0,1]$, with uniformly bounded
dyadic $A_p^d$ or $RH_p^d$~constants, the
\emph{translation-average} $\om(x) := \int_0^1 \om^t(x + t) \,
dt$ is not necessarily a doubling weight. Therefore $\om$ need
not be in true~$A_p$, nor in true~$RH_p$.

The main result of the current paper is that, by contrast, the
\emph{geometric-arithmetic average}~$\Om(x) := \exp \{\int_0^1
\log\om^t(x + t) \, dt\}$ always turns a measurably varying
family of suitably normalized dyadic $A_p^d$ weights into an
$A_p$ weight (for arbitrary $p$ with $1 \leq p \leq \infty$),
and a measurably varying family of suitably normalized dyadic
$RH_p^d$ weights into an $RH_p$ weight (for arbitrary $p$ with
$1 < p \leq \infty$).

\begin{thm}\label{thm:Apexplog}
  Fix $p$ with $1 \le p \le \infty$. Let $\{\om^t\}_{t\in
  [0,1]}$ be a family of dyadic~$A_p$ weights on the circle~$\T$,
  $\om^t\in A_p^d(\T)$, such that
  \begin{enumerate}
  \item[{\rm (i)}] the mapping $t\mapsto\om^t$ is
      measurable,

  \item[{\rm (ii)}] an appropriate average of the
      logarithms of the weights $\om^t$ is finite:
      \[
        \int_0^1 \!\!\! \int _0^1 \big|\log \om^t(x)\big| \, dx \, dt
        < \infty,
      \]
      and

  \item[{\rm (iii)}] the $A_p^d$~constants $A_p^d(\om^t)$
      are uniformly bounded, independent of $t\in [0,1]$.
  \end{enumerate}
  Then the geometric-arithmetic average
  \[
    \Om(x) = \exp\left\{\int_0^1 \log\om^t(x+t) \, dt\right\}
  \]
  of the dyadic weights~$\om^t$ belongs to $A_p$ on~$\T$. Moreover, the $A_p$
  constant of~$\Om$ depends only on~$p$ and on the bound on the
  $A_p^d$~constants of the~$\om^t$.
\end{thm}

\emph{Remark.} There is a simple heuristic motivation for this result. The
weights $\om^t$ are in $A_p^d$, so their logarithms $\log\om^t$
are in~$\bmo_d$. Therefore, as shown in~\cite{GJ} (and
later~\cite{PW} for the one- and two-parameter settings
and~\cite{T} for the general multiparameter setting), the
translation-average~$\log\Om$ of the functions $\log\om^t$ is
in~$\bmo$, and so by the John--Nirenberg Theorem~\cite{JoNi}
sufficiently small powers $\Om^\delta$ of~$\Om$ are in~$A_p$.
It remains to show that $\Om$ itself is in~$A_p$.

\emph{Remark.} Hypothesis (ii) of the theorem is merely a normalization condition:
the fact that each $\om^t(x)$ belongs to $A_p^d$ already implies that
$\log\om^t(x)$ belongs to $L^1(dx)$.

The analogous result to Theorem~\ref{thm:Apexplog} holds for
reverse-H\"older weights.

\begin{thm}\label{thm:RHpexplog}
  Fix $p$ with $1 < p \le \infty$. Let $\{\om^t\}_{t\in [0,1]}$
  be a family of dyadic~$RH_p$ weights on the circle~$\T$, $\om^t\in
  RH_p^d(\T)$, such that hypotheses~\textup{(i)} and~\textup{(ii)}
  of Theorem~\ref{thm:Apexplog} hold, and
  \begin{enumerate}
  \item[{\rm (iii$'$)}] the $RH_p^d$~constants
      $RH_p^d(\om^t)$ are uniformly bounded, independent of
      $t\in [0,1]$.
  \end{enumerate}
  Then the geometric-arithmetic average $\Om(x) =
  \exp\{\int_0^1 \log\om^t(x+t)\, dt\}$ of the dyadic weights $\om^t$
  lies in $RH_p$ on~$\T$. The $RH_p$ constant of $\Om$ depends
  only on~$p$ and on the bound on the $RH_p^d$ constants of the~$\om^t$.
\end{thm}

These results also hold on $\T^k$, with constants that depend on the
dimension~$k$. In addition, they hold in the setting of the polydisc;
see Section~\ref{sec:polydisc}.

We remark that it is not necessary for the integral in~$t$ to
be taken over the whole interval~$[0,1]$. The proofs below go
through without change when the integral is taken over an
arbitrary subset $E \subset [0,1]$ of positive measure.

\section{Definitions and Tools}
\label{sec:definitions}

In this section we collect useful material about doubling
weights, the weight classes~$A_p$ and~$RH_p$, and their
relationship to~$\bmo$. For fuller accounts of the theory of
$A_p$ and $RH_p$ weights, see for example~\cite{GCRF},
\cite{Gar}, \cite{Gra}, and~\cite{CN}.

Let $\T$ denote the unit circle, obtained by identifying the
endpoints of the interval~$[0,1]$. In the definitions of our
averages~$\Om$ and~$\om$, $x + t$ is to be interpreted as $x +
t \mod 1$.

Denote the collection of dyadic subintervals~$I$ of the
circle~$\T$ by $\D = \D[0,1]$:
\[
    \D
    := \{[0,1]\} \cup
        \left\{I = \left[\frac{j}{2^k}, \frac{j+1}{2^k}\right)
        \biggm| k\in\N, j\in\{0, 1, \ldots\,, 2^k - 1\}\right\}.
\]
Throughout the paper, $Q$~denotes a general subinterval
of~$\T$, while $I$, $J$, $K$ and $L$ denote dyadic subintervals
of~$\T$.

The functions we consider are real-valued.

We use the symbol $|E|$ to denote the Lebesgue measure of a
set~$E$, the symbol $\intav_E$ to denote $\frac{1}{|E|}
\int_E$, and the symbol~$f_E$ for the average value~$\intav_E
f$ of a function~$f$ on a set~$E$. The notation $E\subset F$
includes the possibility $E = F$.

\begin{defn}
  \label{def:doubling}
  Let $\omega(x)$ be a nonnegative locally integrable function on
  the circle~$\T$. We say $\om$ is a \emph{doubling weight}
  with \emph{doubling constant}~$C$ if for all intervals
  $Q\subset\T$
  \[
    \int_{\widetilde{Q}} \om(x) \, dx
    \le C \int_Q \om(x) \, dx,
  \]
  where $\widetilde{Q}$ is the \emph{double} of $Q$; that is,
  $\widetilde{Q}$ is the interval with the same midpoint as~$Q$
  and twice the length of~$Q$. We say $\om$ is a \emph{dyadic
  doubling weight} with \emph{dyadic doubling constant}~$C$ if
  the analogous inequality holds for all dyadic intervals
  $I\subset\T$,
  where $\widetilde{I}$ is the \emph{dyadic double} of~$I$:
  that is, $\widetilde{I}$ is unique dyadic interval
  of length~$|\widetilde{I}| = 2|I|$ that contains~$I$.
\end{defn}

The $A_p$ weights were identified by Muckenhoupt~\cite{M} as
the weights~$\om$ for which the Hardy--Littlewood maximal
function is bounded from~$L^p(d\mu)$ to itself, where $d\mu =
\om(x) \, dx$. Here we give the definitions of the classes
$A_p$ and $RH_p$ on the circle~$\T$; the equivalent definitions
hold on~$\R$ and on (one-parameter)~$\R^n$. We delay the
corresponding definitions for the polydisc setting until
Section~\ref{sec:polydisc}.

\begin{defn}
  \label{def:Ap}
  Let $\omega(x)$ be a nonnegative locally integrable function
  on the circle~$\T$. For real~$p$ with $1 < p < \infty$, we
  say $\omega$ is an $A_p$ \emph{weight}, written $\omega\in
  A_p$, if
  \[
    A_p(\om)
    := \sup_Q \left(\intav_Q \omega\right)
    \left(\intav_Q
      \left(\dfrac{1}{\omega}\right)^{1/(p-1)}\right)^{p-1}
    < \infty.
  \]
  For $p = 1$, we say $\omega$ is an $A_1$ \emph{weight},
  written $\omega\in A_1$, if
  \[
    A_1(\om)
    := \sup_Q \left(\intav_Q \omega\right)
    \left(\dfrac{1}{\essinf_{x\in Q} \omega (x)}\right)
    < \infty.
  \]
  For $p = \infty$, we say $\omega$ is an \emph{$A$-infinity
  weight}, written $\omega\in A_\infty$, if
  \[
    A_\infty(\om)
    := \sup_Q \left(\intav_Q \omega\right)
    \exp\left(\intav_Q \log \left(\frac{1}{\omega}\right) \right)
    < \infty.
  \]
  Here the suprema are taken over all intervals~$Q\subset\T$.
  The quantity $A_p(\om)$ is called the \emph{$A_p$~constant
  of~$\om$}.

  The \emph{dyadic $A_p$ classes} $A_p^d$ for $1 \le p \le
  \infty$ are defined analogously, with the
  suprema~$A_p^d(\om)$ being taken over only the dyadic
  intervals~$I\subset\T$.
\end{defn}

\begin{defn}
  \label{def:RHp}
  Let $\omega(x)$ be a nonnegative locally integrable function
  on the circle~$\T$. For real~$p$ with $1 < p < \infty$, we
  say $\omega$ is a \emph{reverse-H\"older-$p$ weight}, written
  $\omega\in RH_p$ or $\om\in B_p$, if
  \[
    RH_p(\om)
    := \sup_Q \left(\intav_Q \omega^p\right)^{1/p}
    \left(\intav_Q \omega\right)^{-1}
    < \infty.
  \]
  For $p = \infty$, we say $\omega$ is a
  \emph{reverse-H\"older-infinity weight}, written $\omega\in
  RH_\infty$ or $\om\in B_\infty$, if
  \[
    RH_\infty(\om)
    := \sup_Q \left(\esssup_{x\in Q} \om\right)
    \left(\intav_Q \omega\right)^{-1}
    < \infty.
  \]
  Here the suprema are taken over all intervals~$Q\subset\T$.
  The quantity $RH_p(\om)$ is called the \emph{$RH_p$~constant}
  of~$\om$.

  For $1 < p \le \infty$, we say $\omega$ is a \emph{dyadic
  reverse-H\"older-$p$ weight}, written $\om\in RH_p^d$ or $\om
  \in B_p^d$, if the analogous condition $\sup_{I\in\D}
  \left(\intav_I \omega^p\right)^{1/p} \left(\intav_I
  \omega\right)^{-1} < \infty$ or $\sup_{I\in\D}
  \left(\esssup_{x\in I} \om\right) \left(\intav_I
  \omega\right)^{-1} < \infty$ holds with the supremum being
  taken over only the dyadic intervals $I\subset\T$, and if in
  addition $\om$ is a dyadic doubling weight. We define the
  \emph{$RH_p^d$~constant} $RH_p^d(\om)$~of~$\om$ to be the
  larger of this dyadic supremum and the dyadic doubling
  constant.
\end{defn}

The $A_p$ inequality (or the $RH_p$ inequality) implies that
the weight~$\om$ is doubling, and the dyadic $A_p$ inequality
implies that~$\om$ is dyadic doubling. However, the dyadic
$RH_p$ inequality does not imply that $\om$ is dyadic doubling,
which is why the dyadic doubling assumption is needed in the
definition of~$RH_p^d$.

The $A_p$ classes are nested and increasing with~$p$, while the
$RH_p$~classes are nested and decreasing with~$p$. Moreover,
\[
    A_q(\om) \le A_p(\om), \qquad
    \text{for $1 \le p < q \le \infty$,}
\]
and
\[
    RH_p(\om) \le RH_q(\om), \qquad
    \text{for $1 < p < q \le \infty$.}
\]
Also
\[
    A_\infty
    = \bigcup_{p \ge 1} A_p
    = \bigcup_{q > 1} RH_q, \qquad
    A_1
    \subsetneq \bigcap_{p > 1} A_p, \qquad
    RH_\infty
    \subsetneq \bigcap_{p > 1} RH_p.
\]
The dyadic versions of the assertions in this paragraph also
hold.

The example $w(x) = \left(\log (1/|x|)\right)^{-1}$ (for $x$ near zero)
cited in~\cite{JoNe} shows that $A_1$ is a proper subset of $\cap_{p > 1} A_p$.
The example $\om(x) = \max\{\log(1/|x|),1\}$ given in~\cite{CN} shows that
$RH_\infty$ is a proper subset of $\cap_{p > 1} RH_p$. However, as
noted in~\cite{CN}, if a weight~$\om$ is in $A_p$ for each $p
> 1$ and if the constants~$A_p(\om)$ are uniformly bounded,
then $\om\in A_1$; and the corresponding statement holds for~$RH_p$
and~$RH_\infty$.

As noted above, for a nonnegative locally integrable
function~$\om$,
\begin{align*}
    \om~\text{is in}~A_\infty
    & \Longleftrightarrow
    \om~\text{is in}~A_p~\text{for some}~p \in [1,\infty) \\
    & \Longleftrightarrow
    \om~\text{is in}~RH_q~\text{for some}~q \in (1,\infty).
\end{align*}
In the first equivalence the $A_\infty$ constant depends only
on the $A_p$ constant and on~$p$, which in turn depend only on
the $A_\infty$~constant. Similarly, the $A_\infty$ constant
depends only on the $RH_q$ constant and on~$q$, which depend
only on the $A_\infty$~constant. See for
example~\cite[Theorem~9.3.3]{Gra} where the constants in these
and other characterizations of~$A_\infty$ are carefully
analyzed. The analogous statements hold for the dyadic
classes~$A_\infty^d$, $A_p^d$, and $RH_p^d$.


The classes of $A_p$ and $RH_p$ weights can be characterized by
conditions on the oscillation of the logarithm of the weight,
as follows.

\begin{lem}\label{lem:ApRHpcharn}
  Let $\omega$ be a nonnegative locally integrable function
  on~$\T$. Let $\p := \log \omega$. Then the following five
  statements hold.
  \begin{enumerate}
  \item[\textup{(a)}] $\om$ is in $A_\infty$ if and only if
  \begin{equation}\label{eqn:Ainfinity}
    \sup_Q \intav_Q \exp\{\p(x) - \p_Q\} \, dx
    < \infty.
  \end{equation}

  \item[\textup{(b)}] For $1 < p < \infty$, $\om$ is in
      $A_p$ if and only if inequality~\eqref{eqn:Ainfinity}
      holds and also
  \begin{equation}\label{eqn:Ap}
    \sup_Q \intav_Q
    \exp\left\{\frac{-(\p(x) - \p_Q)}{p - 1}\right\} \, dx
    < \infty.
  \end{equation}

  \item[\textup{(c)}] $\om$ is in $A_1$ if and only if
      inequality~\eqref{eqn:Ainfinity} holds and also
  \begin{equation}\label{eqn:A1}
    \sup_Q \big[\p_Q - \essinf_{x\in Q} \p(x)\big]
    < \infty.
  \end{equation}

  \item[\textup{(d)}] $\om$ is in $RH_\infty$ if and only
      if
  \begin{equation}\label{eqn:RHinfinity}
    \sup_Q \big[\esssup_{x\in Q} \p(x) - \p_Q\big] <
    \infty.
  \end{equation}

  \item[\textup{(e)}] For $1 < p < \infty$, $\om$ is in
      $RH_p$ if and only if
  \begin{equation}\label{eqn:RHp}
    \sup_Q \intav_Q \exp\{p(\p(x) - \p_Q)\} \, dx
    < \infty.
  \end{equation}
  \end{enumerate}

  In each part, the value of $A_p(\om)$ or $RH_p(\om)$ depends
  on the value(s) of the supremum (suprema) in the
  characterization given and, when $1 < p < \infty$, also
  on~$p$. Conversely, the value(s) of the supremum (suprema)
  depend on the value of $A_p(\om)$ or $RH_p(\om)$ and, when $1
  < p < \infty$, also on~$p$.

  Taking the suprema over only dyadic intervals~$I\subset\T$,
  the dyadic analogues of parts~(a)--(e) hold for the dyadic
  classes $A_\infty^d$, $A_p^d$, $A_1^d$, $RH_\infty^d$,
  and~$RH_p^d$, except that in parts~(d) and~(e) one needs in
  addition to inequality~\eqref{eqn:RHinfinity}
  or~\eqref{eqn:RHp} the extra hypothesis that $\om$ is dyadic
  doubling. The dependence of the constants in the dyadic
  case is the same as in the continuous case.
\end{lem}

Parts~(b) and~(c) appear in~\cite{Gar}, \cite{Gra},
and~\cite{GCRF}, for example, and part~(d) is in~\cite[Cor~4.6]{CN}.
Inequality~\eqref{eqn:A1} says that $\om$ is in $A_1$ when its logarithm~$\p$
belongs to the space $\blo$ of functions of \emph{bounded lower oscillation},
while inequality~\eqref{eqn:RHinfinity} says that $\om$ is in~$RH_\infty$ when
$-\p$ belongs to $\blo$.

\begin{proof} Let $C_1$, $C_2$, $C_3$, $C_4$, and $C_5$ be the
suprema in inequalities~\eqref{eqn:Ainfinity}, \eqref{eqn:Ap},
\eqref{eqn:A1}, \eqref{eqn:RHinfinity}, and~\eqref{eqn:RHp}
respectively.

(a) It is immediate that $A_\infty(\om) = C_1$, since for each
interval~$Q$ the $A_\infty$ quantity is
\begin{align*}
    \left(\intav_Q \om(x) \, dx\right)
        \exp\left\{\intav_Q \log \frac{1}{w(x)} \, dx\right\}
    & = \intav_Q \exp\{\p(x) - \p_Q\} \, dx.
\end{align*}

(b) We show that $C_1 \le A_p(\om)$, $C_2 \le
A_p(\om)^{1/(p-1)}$, and $A_p(\om) \le C_1 C_2^{p-1}$. Let
\[
    \psi
    := \log \left[\left(\frac{1}{\om}\right)^{1/(p-1)}\right]
    = \frac{-\p}{p-1}.
\]
Let $Q$ be an interval in~$\T$. Then by Jensen's inequality,
  \begin{align*}
    \intav_Q \exp\{\p(x) - \p_Q\} \, dx
    & = \left(\intav_Q w\right) \exp\left((p-1) \intav_Q \psi(x) \, dx\right) \\
    & \le \left(\intav_Q w\right) \left[\intav_Q \exp \psi(x) \,
         dx\right]^{p-1} \\
    & = \left(\intav_Q w\right)
        \left[\intav_Q \left(\frac{1}{\om}\right)^{1/(p-1)}\right]^{p-1} \\
    & \le A_p(\om).
  \end{align*}
  Thus inequality~\eqref{eqn:Ainfinity} holds with $C_1 = A_p(\om)$.
  Similarly, by Jensen's inequality,
  \begin{align*}
    \intav_Q \exp\left\{\frac{-(\p(x) - \p_Q)}{p-1}\right\} \, dx
    & = \bigg(\intav_Q \left(\frac{1}{\om}\right)^{1/(p-1)}\bigg)
        \bigg[\exp\left(\intav_Q \p \right)\bigg]^{1/(p-1)} \\
    & \leq \left[\bigg(\intav_Q \left(\frac{1}{\om}\right)^{1/(p-1)}
        \bigg)^{p-1}
        \left(\intav_Q \om \right)\right]^{1/(p-1)} \\
    & \le A_p(\om)^{1/(p-1)}.
  \end{align*}
  Thus inequality~\eqref{eqn:A1} holds with $C_2 =
  A_p(\om)^{1/(p-1)}$.

  For the converse,
  \begin{align*}
    \left(\intav_Q \om \right)
        \bigg(\intav_Q \left(\frac{1}{\om}\right)^{1/(p-1)}\bigg)^{p-1}
    & = \left(\intav_Q \exp\p(x) \, dx\right)
        \bigg(\intav_Q \exp\psi(x) \, dx\bigg)^{p-1}
        e^{-\p_Q} e^{-(p-1)\psi_Q}\\
    & = \left(\intav_Q \exp\{\p(x) - \p_Q\} \, dx\right)
        \bigg(\intav_Q \exp\{\psi(x) - \psi_Q\} \, dx\bigg)^{p-1} \\
    & \le C_1 C_2^{p-1},
  \end{align*}
  and thus $A_p(\om) \le C_1 C_2^{p-1}$.

(c) We show that $C_1 \le A_1(\om)$, $C_3 \le \log A_1(\om)$,
and $A_1(\om) \le C_1 \exp C_3$. If $\om$ is in~$A_1$, then for
each interval~$Q$ in~$\T$ we have
\[
  \intav_Q e^{\p(x)} \, dx
  = \intav_Q w(x) \, dx
  \leq A_1(\om) \essinf_{x\in Q} w(x)
  \leq A_1(\om) e^{\p_Q}.
\]
It follows that
\[
  \intav_Q e^{\p(x) - \p_Q} \, dx
  \leq A_1(\om).
\]
Thus $\p$ satisfies inequality~\eqref{eqn:Ainfinity} with
constant~$C_1 \le A_1(\om)$.

By Jensen's inequality and the $A_1$~property,
\begin{equation}
  e^{\p_Q}
  \leq \intav_Q w(x) \, dx
   \leq A_1(\om) \exp\big\{\essinf_{x\in Q} \p(x)\big\}.
\end{equation}
Therefore
\[
  \p_Q
  \leq \log A_1(\om) + \essinf_{x\in Q} \p(x).
\]
Thus $\p$ satisfies inequality~\eqref{eqn:A1} with
constant~$C_3 = \log A_1(\om)$.

Now suppose that $\p$ satisfies
inequalities~\eqref{eqn:Ainfinity} and~\eqref{eqn:A1}. Then for
each interval~$Q$,
\begin{align*}
  \intav_Q w(x) \, dx
  &= \intav_Q e^{\p(x)} \, dx
  \leq C_1 \, e^{\p_Q} \\
  &\leq C_1 \exp\big\{C_3 + \essinf_{x\in Q} \p(x)\big\}
  = C_1 e^{C_3} \essinf_{x\in Q} w(x).
\end{align*}
Thus $\om$ satisfies the $A_1$~property with constant $A_1(\om)
\le C_1 e^{C_3}$.

(d) We show that $C_4 \le \log (RH_\infty(\om) A_\infty(\om))$
and $RH_\infty(\om) \le e^{C_4}$, and that the bound on~$C_4$
depends only on~$RH_\infty(\om)$.

Suppose $\om$ is in $RH_\infty$. Then $\om$ is in~$A_\infty$, so
inequality~\eqref{eqn:Ainfinity} holds with $C_1 = A_\infty(\om)$.
Further, $\om$ is in every~$RH_p$ for $p\in (1,\infty)$, and
$A_\infty(\om)$ depends only on $RH_p(\om)$, while $RH_p(\om) \le
RH_\infty(\om)$. Thus $A_\infty(\om)$ depends only
on~$RH_\infty(\om)$. Now for each interval~$Q$ in~$\T$,
inequality~\eqref{eqn:Ainfinity} implies that
\[
    \esssup_{x\in Q} \om(x)
    \le RH_\infty(\om) \intav_Q e^{\p(x)} \, dx
    \le RH_\infty(\om) A_\infty(\om) e^{\p_Q}.
\]
Taking logarithms, we see that
\[
    \esssup_{x\in Q} \p(x)
    \le \p_Q + \log(RH_\infty(\om) A_\infty(\om)),
\]
and so inequality~\eqref{eqn:RHinfinity} holds with $C_4 =
\log(RH_\infty(\om) A_\infty(\om))$.

Conversely, if inequality~\eqref{eqn:RHinfinity} holds, then by
Jensen's inequality
\[
    \esssup_{x\in Q} \om(x)
    \le \exp\{C_4 + \p_Q\}
    \le e^{C_4} \intav_Q \exp{\p(x)} \, dx
    = e^{C_4} \intav_Q \om.
\]
Thus $RH_\infty(\om) \le e^{C_4}$.

(e) We show that $C_5 \le RH_p(\om) A_\infty(\om)$ and
$RH_p(\om) \le C_5^{1/p}$, and that $C_5$ depends only on~$p$
and on~$RH_p(\om)$. In terms of~$\p$, the $RH_p$ expression for
a given interval~$Q$ is
\[
    \left(\intav_Q \om^p\right)^{1/p}
        \left(\intav_Q \om\right)^{-1}
    = \left(\intav_Q \exp\{p(\p(x) - \p_Q)\} \, dx\right)^{1/p}
        \left(\intav_Q \exp\{\p(x) - \p_Q\} \, dx\right)^{-1}.
\]
Also, if $\om$ is in $RH_p$, then $\om$ is in~$A_\infty$ and
$A_\infty(\om)$ depends only on~$p$ and on~$RH_p(\om)$. It
follows that
\[
    \left(\intav_Q \exp\{p(\p(x) - \p_Q)\} \, dx\right)^{1/p}
    \le RH_p(\om) \intav_Q \exp\{\p(x) - \p_Q\} \, dx
    \le RH_p(\om) A_\infty(\om).
\]
Thus inequality~\eqref{eqn:RHp} holds with $C_5 \le RH_p(\om)
A_\infty(\om)$, and~this bound depends only on~$p$ and
on~$RH_p(\om)$.

By Jensen's inequality, $\intav_Q \exp\{\p(x) - \p_Q\} \, dx
\geq 1$. Thus if inequality~\eqref{eqn:RHp} holds, then
\[
    \left(\intav_Q \exp\{p(\p(x) - \p_Q)\} \, dx\right)^{1/p}
        \left(\intav_Q \exp\{\p(x) - \p_Q\} \, dx\right)^{-1}
    \le C_5^{1/p}.
\]
Thus $\om$ is in $RH_p$ and $RH_p(\om) \le C_5^{1/p}$.

The same arguments go through for the dyadic classes~$A_p^d$
and~$RH_p^d$.
\end{proof}

Muckenhoupt's $A_p$ weights are closely related to functions of
bounded mean oscillation.

\begin{defn}
  A real-valued function $f\in L^1(\T)$ is in the space
  $\bmo(\T)$ of functions of \emph{bounded mean oscillation} on
  the circle if its $\bmo$ norm is finite:
  \[
    \Vert f \Vert_* := \sup_{Q\subset\T} \,
    \intav_Q \, |f(x) - f_Q| \, dx < \infty.
  \]
   \emph{Dyadic $\bmo$} of the circle, written $\bmo_d(\T)$, is
  the space of functions that satisfy the corresponding
  estimate where the supremum is taken over all dyadic
  subintervals $I\in\D$ of $[0,1]$.  The dyadic $\bmo$ norm of
  $f$ is denoted $\Vert f \Vert_d$.
\end{defn}

Elements of $\bmo$ (or $\bmo_d$) that differ only by an
additive constant are equivalent; thus $\bmo$ and $\bmo_d$ are
subspaces of~$L^1_\text{loc}/\R$.

For $1 \le p \le \infty$, if $\om$ is in~$A_p$ then $\p :=
\log\om$ is in $\bmo$, with $\bmo$ norm depending only on the
$A_p$ constant $A_p(\om)$. See for example~\cite[p.409]{GCRF}.
The same is true for $RH_p$ weights. Specifically, we have the
following result; we omit the proof.





\begin{lem}\label{lem:AptoBMO}
    Suppose $1 \le p \le \infty$. If $\om$ is in $A_p$ then $\p :=
    \log\om$ is in~$\bmo$. For $1 < p < \infty$,
    \[
        ||\p||_*
        \le A_p(\om) +(p-1)A_p(\om)^{1/(p-1)}.
    \]
    For $p=1$, $||\p||_* \le 2A_1(\om)$.
    For $p = \infty$, $||\p||_*$ depends only on~$A_\infty(\om)$.
    For $1 < p \le \infty$, if $\om$ is in $RH_p$ then $\p :=
    \log\om$ is in~$\bmo$, with $||\p||_*$ depending only on
    $RH_{p}(\om)$ and on~$p$. The analogous statements hold in
    the dyadic setting.
\end{lem}


We use a characterization of the dyadic $\bmo$ functions on the
circle in terms of the size of Haar coefficients. The Haar
function~$h_I$ associated with the dyadic interval $I$ is given
by $h_I(x) = |I|^{-1/2}$ if $x$ is in the left half of~$I$,
$h_I(x) = -|I|^{-1/2}$ if $x$ is in the right half of~$I$, and
$h_I = 0$ otherwise.
The Haar coefficient over $I$ of $f$ is $f_I = (f,h_I) := \int_I \, f(x) h_I(x) \, dx$.
The Haar series for $f$ is
\[
  f(x) := \sum_{I\in\D} \, (f,h_I) \, h_I(x),
\]
and the $L^2$-norm of $f$ is given in terms of the Haar
coefficients by
\[
  \Vert f \Vert_2
  = \bigg[\sum_{J\in\D} \, (f,h_J)^2\bigg]^{1/2}.
\]

It follows from the John--Nirenberg Theorem~\cite{JoNi} that
for $p \geq 1$, for $f$ in~$L^1(\T)$ the expression
\[
  \Vert f \Vert_{d,p}
  := \sup_{I\in\D} \,
    \left(\intav_I \, |f(x) - (f)_I|^p \, dx\right)^{1/p}
\]
is comparable to the dyadic $\bmo$ norm $\Vert f \Vert_d$.

In particular, a function $f\in L^1(\T)$ of mean value zero is
in $\bmo_d(\T)$ if and only if there is a constant $C$ such
that for all $I\in\D$,
\begin{equation}\label{eqn:haarcarleson}
  \sum_{J\subset I, J\in\D} \, (f,h_J)^2 \le C|I|.
\end{equation}
Moreover, the smallest such constant $C$ is equal to $\Vert f
\Vert_{d,2}^2$.

Since the sum in inequality~\eqref{eqn:haarcarleson} ranges
over dyadic intervals only, there is no need to restrict the
interval~$I$ itself to be dyadic.

\section{Proofs of Theorems~\ref{thm:Apexplog} and~\ref{thm:RHpexplog}}
\label{sec:explogApRHp}

We begin this section with three lemmas, which we then use to
prove the geometric-arithmetic averaging result for both~$A_p$ and $RH_p$.

Lemma \ref{lem:CACBestimates} below gives an estimate on Haar
expansions of $\bmo_d$ functions.
Lemmas~\ref{lem:betainequality}
and~\ref{lem:A1RHinfinityinequality}, which rely on the
estimates~\eqref{eqn:ca} and~\eqref{eqn:cb} from
Lemma~\ref{lem:CACBestimates}, will allow us to pass from the
dyadic versions to the non-dyadic versions of the inequalities
that characterize $A_p$ and $RH_p$.

Throughout this section we use the following notation.
Let $\D_n := \{I\in\D \bigm| |I| = 2^{-n}\}$ be the collection of dyadic
intervals of length~$2^{-n}$, for $n = 0$, 1, 2, \dots\  Expanding each $\p^t$
in Haar series, we have
\[ \p(x) =
  \int_0^1 \sum_{J\in \D}(\p^t,h_J)h_J(x+t)\,dt
  = \sum_{n=0}^{\infty}\int_0^1 \,
      \sum_{J\in\D_n} \, (\pa,h_J) \, h_J(x + t) \, dt
  = \sum_{n=0}^\infty \p_n(x),
\]
so that $\p_n$ is the translation-average over~$t$ of the
slices at scale~$2^{-n}$ of the Haar expansions for the
functions~$\p^t$.

Fix an interval $Q\subset\T$; this~$Q$ need not necessarily be
dyadic.
Split the sum for $\p(x)$, at the scale of $|Q|$, into two
parts $\p_A$ and $\p_B$ in which the dyadic intervals~$J$ are
respectively small and large compared with~$Q$:
\[
  \p = \p_A + \p_B,
  \qquad \p_A(x) := \sumA \, \p_n(x),
  \qquad \p_B(x) := \sumB \, \p_n(x).
\]

The following result is proved in the course of the proof of
Theorem~2 of~\cite{PW}.

\begin{lem}\label{lem:CACBestimates}
Suppose that $\{\p^t\}_{t\in[0,1]}$ is a family of dyadic \bmo\
functions on~$\T$, $\p^t\in\bmo_d$, such that
\begin{enumerate}
  \item[{\rm (i)}] the mapping $t\mapsto\p^t$ is
      measurable,

  \item[{\rm (ii)}] the $\bmo_d$~constants $||\p^t||_d$ are
      uniformly bounded, independent of $t\in [0,1]$, and

  \item[{\rm (iii)}] for each $t\in[0,1]$, the function
      $\p^t$ has mean value zero on~$\T$.
  \end{enumerate}
Then there are constants $C_A$ and $C_B$ depending on the bound
on the $\bmo_d$ constants~$\|\p^t\|_d$, and independent of~$Q$,
such that for each interval~$Q\subset\T$ and for each
point~$x_0\in Q$,
\begin{align}
  \intq |\p_A(x)|^2 \, dx
  & \le C_A, \label{eqn:ca} \\
  \intq |\p_B(x) - \p_B(x_0)| \, dx
  & \le C_B. \label{eqn:cb}
\end{align}
\end{lem}



\begin{lem}\label{lem:betainequality}
  Let $\beta$ be a real number. Suppose $\{\om^t\}_{t\in
  [0,1]}$ is a family of nonnegative locally integrable
  functions on~$\T$ such that for all $\p^t:=\log\om^t$,
  hypotheses \textup{(i)--(iii)} of
  Lemma~\ref{lem:CACBestimates} hold. Let $\p(x)
    = \log\Om(x)
    := \int_0^1 \log\omega^t(x+t)\, dt$.
  Suppose there is a constant~$C^d(\beta)$ such that for all
  $t\in[0,1]$ and for all dyadic intervals~$I \subset \T$,
  \begin{equation}\label{eqn:betainequalitydyadic}
    \intav_I \exp\big[\beta(\p^t(x) - \p^t_I)\big] \, dx
    \le C^d(\beta).
  \end{equation}
  Then there is a constant~$C(\beta)$, depending only on
  $C^d(\beta)$, such that
  \begin{equation}\label{eqn:betainequality}
    \intav_Q \exp\big[\beta(\p(x) - \p_Q)\big] \, dx
    \le C(\beta)
  \end{equation}
  for all intervals~$Q\subset\T$.
\end{lem}

In fact, for many choices of $\beta$, hypothesis (ii) of
Lemma~\ref{lem:CACBestimates} is implied by
inequality~\eqref{eqn:betainequalitydyadic}, together with
Lemmas~\ref{lem:ApRHpcharn} and~\ref{lem:AptoBMO}. This is made
clear in the proof of Theorems~\ref{thm:Apexplog}
and~\ref{thm:RHpexplog}.

\begin{proof}
We first establish an inequality controlling the exponentials
of the Haar expansions of the $\p^t$. By
inequality~\eqref{eqn:betainequalitydyadic}, for each dyadic
interval~$I\subset\T$ we have
\begin{align}\label{eqn:expCarleson}
  \intav_I \exp\Big[\beta\sum_{J\subset I, J\in\D}(\p^t,h_J)h_J(x)\Big]\, dx
  & = \intav_I \exp\Big[\beta\sum_{J\in\D}(\p^t,h_J)h_J(x)
    - \beta\sum_{J\varsupsetneq I, J\in\D}(\p^t,h_J)h_J(x)\Big]\, dx \nonumber\\
  & = \intav_I \exp\big[\beta(\p^t(x) - \p^t_I)\big] \, dx \nonumber\\
  & \leq C^d(\beta).
\end{align}
We have used the elementary fact that the average~$f_I$ of a
function~$f \in L^1(\T)$ over an interval~$I$ containing the
point~$x$ can be written as
\[
   f_I
   = \intav_I \sum_{J\in\D} (f,h_J) h_J(s) \, ds
   = \sum_{J\supsetneq I, J\in\D} (f,h_J) h_J(x).
\]

Fix an interval~$Q\subset\T$, not necessarily dyadic. For
each~$x\in Q$ we have
\begin{align*}
  \left|\p_B(x) - \intav_Q \p(s) \, ds \right|
  & \le \intav_Q |\p_B(x) - \p(s)| \, ds \\
  & \le \intav_Q |\p_A(s)| \, ds + \intav_Q |\p_B(x) - \p_B(s)| \, ds \\
  & \le \sqrt{C_A} + C_B,
\end{align*}
by Cauchy--Schwarz and the estimates~\eqref{eqn:ca}
and~\eqref{eqn:cb} from Lemma~\ref{lem:CACBestimates}.
Therefore
\begin{align*}
  \intav_Q \exp\big[\beta(\p(x) - \p_Q)\big] \, dx
  & = \intav_Q \exp\big[\beta\p_A(x)\big] \,
    \exp\big[\beta(\p_B(x) - \p_Q)\big] \, dx \\
  & \le \exp\big[\beta(\sqrt{C_A} + C_B)\big] \,
    \intav_Q \exp\big[\beta\p_A(x)\big] \, dx.
\end{align*}
Thus it suffices to bound the quantity
\begin{align*}
  \intav_Q \exp[\beta\p_A(x)] \, dx
  & = \intav_Q \exp \bigg[\beta\sum_{2^{-n} < |Q|} \int_0^1 \sum_{J\in\D_n}
    (\p^t,h_J) h_J(x + t) \, dt\bigg] \, dx \\
  & = \intav_Q \exp \bigg[\int_0^1 \beta \sum_{|J| < |Q|, J\in\D}
    (\p^t,h_J) h_J(x + t) \, dt\bigg] \, dx \\
  & \le \int_0^1 \!\!\! \intav_Q \exp \bigg[\beta\sum_{|J| < |Q|, J\in\D}
    (\p^t,h_J) h_J(x + t)\bigg] \, dx \, dt \\
  & = \int_0^1 \!\!\! \intav_{Q+t} \exp \bigg[\beta\sum_{|J| < |Q|, J\in\D}
    (\p^t,h_J) h_J(x)\bigg] \, dx \, dt.
\end{align*}
We have used Jensen's inequality and Tonelli's Theorem.

In order to apply inequality~\eqref{eqn:expCarleson}, we want
to replace the interval~$Q$ in the preceding expression by
appropriate dyadic intervals. Fix~$t$. There are two adjacent
dyadic intervals $K$ and~$L$ such that $Q + t \subset K \cup L$
and $|Q| \le |K| = |L| < 2|Q|.$ Then
\begin{align*}
  \intav_Q \exp[\beta\p_A(x)] \, dx
  & \le \int_0^1 \frac{1}{|Q|} \int_{Q+t}
    \exp \bigg[\beta\sum_{|J| < |Q|,
    J\in\D} (\p^t,h_J) h_J(x)\bigg] \, dx \, dt \\
  & \le \int_0^1 \frac{2}{|K|} \int_{K\cup L}
    \exp \bigg[\beta\sum_{|J| < |K|,
    J\in\D} (\p^t,h_J) h_J(x)\bigg] \, dx \, dt \\
  & = 2\int_0^1 \bigg\{\frac{1}{|K|} \int_{K}
    \exp \bigg[\beta\sum_{|J| < |K|,
    J\in\D} (\p^t,h_J) h_J(x)\bigg] \, dx \\
  &   {} \qquad\qquad + \frac{1}{|L|} \int_{L}
    \exp \bigg[\beta\sum_{|J| < |L|,
    J\in\D} (\p^t,h_J) h_J(x)\bigg] \, dx\bigg\} \, dt \\
  & \le 2 \int_0^1 C^d(\beta) + C^d(\beta) \, dt
  = 4 \, C^d(\beta),
\end{align*}
by inequality~\eqref{eqn:expCarleson}. Thus
\[
  \intav_Q \exp\big[\beta(\p(x) - \p_Q)\big] \, dx
  \le 4 \, C^d(\beta) \exp\left[\beta\left(\sqrt{C_A} + C_B\right)\right].
\]
Taking the supremum over all intervals~$Q\subset\T$, we see
that inequality~\eqref{eqn:betainequality} holds for~$\p =
\log\Om$, with constant~$C(\beta) = 4 \, C^d(\beta)
\exp\left[\beta\left(\sqrt{C_A} + C_B\right)\right]$.
\end{proof}

\begin{lem}\label{lem:A1RHinfinityinequality}
  Suppose $\{\om^t\}_{t\in [0,1]}$ is a family of nonnegative
  locally integrable functions on~$\T$ such that for all
  $\p^t:=\log\om^t$, hypotheses \textup{(i)} and \textup{(iii)}
  of Lemma \ref{lem:CACBestimates} hold. As before, let $\p(x)
    = \log\Om(x)
    := \int_0^1 \log\omega^t(x+t)\, dt$.
  \begin{enumerate}
  \item[(a)] Suppose there is a constant $C_3^d$ such that
  for all $t \in [0,1]$ and for all dyadic intervals $I \subset \T$,
  \[
    \big[\p^t_I - \essinf_{x\in I} \p^t(x)\big]
    \le C_3^d.
  \]
  Then there is a constant
  $C_3$ depending only on $C_3^d$ such that
  \[
    \big[\p_Q - \essinf_{x\in Q} \p(x)\big]
    \le C_3
  \]
  for all intervals~$Q\subset\T$.
  \item[(b)] Similarly, if there is a constant $C_4^d$
  such that
  for all $t \in [0,1]$ and for all dyadic intervals $I \subset \T$,
  \[
    \big[\esssup_{x\in I} \p^t(x) - \p^t_I\big]
    \le C_4^d,
  \]
  then there is a constant $C_4$ depending only on $C_4^d$
  such that
  \[
    \big[\esssup_{x\in Q} \p(x) - \p_Q\big]
    \le C_4
  \]
  for all intervals~$Q\subset\T$.
  \end{enumerate}
\end{lem}

\begin{proof}
Observe that hypothesis (ii) of Lemma \ref{lem:CACBestimates}
follows from the assumption in part~(a) or the assumption in
part~(b), together with Lemmas~\ref{lem:ApRHpcharn}
and~\ref{lem:AptoBMO}. In particular, the $\bmo_d$~constants
$\|\p^t\|_d$ depend only on $C_3^d$ or~$C_4^d$.

(a) For each dyadic interval $I\subset\T$ and for a.e.~$x\in
I$,
\begin{equation}\label{eqn:bloHaar}
  C_3^d
  \ge \p^t_I - \p^t(x)
  = - \sum_{J\subset I, J\in\D} (\p^t,h_J) h_J(x).
\end{equation}
Fix an interval~$Q\subset\T$, not necessarily dyadic. For $x\in
Q$ consider the quantity
\begin{align*}
  \p_Q - \p(x)
  & = \intav_Q \big[\p_A(s) + \p_B(s)\big] \, ds
        - \big[\p_A(x) + \p_B(x)\big] \\
  & = \underbrace{\intav_Q \p_A(s) \, ds}_{(\text{I})}
        + \underbrace{\left[\intav_Q \p_B(s) \, ds
          - \p_B(x)\right]}_{(\text{II})}
        + \underbrace{\bigg[-\p_A(x)\bigg]}_{(\text{III})}.
\end{align*}
We bound the terms $(\text{I})$, $(\text{II})$, and $(\text{III})$
separately. By Lemma~\ref{lem:CACBestimates},
\[
  (\text{I})
  \le \left[\intav_Q |\p_A(s)|^2 \, ds\right]^{1/2}
  \le \sqrt{C_A}
\]
and
\[
  (\text{II})
  \le \intav_Q |\p_B(s) - \p_B(x)| \, ds
  \le C_B
\]
for all~$x\in Q$. Also, for $x\in Q$ and $t\in[0,1]$ there is a
unique dyadic interval~$I_{x,t}$ such that $x + t \in I_{x,t}$
and $|Q|/2 \le |I_{x,t}| < |Q|$. Then
\begin{align*}
  {}- \p_A(x)
  & = - \sumA \int_0^1 \sum_{J\in\D_n}
        (\p^t,h_J) h_J(x + t) \, dt \\
  & = - \int_0^1 \sum_{|J| < |Q|, J\in\D}
        (\p^t,h_J) h_J(x + t) \, dt \\
  & = - \int_0^1 \sum_{J\subset I_{x,t}, J\in\D}
        (\p^t,h_J) h_J(x + t) \, dt \\
  & \le C_3^d
\end{align*}
for a.e.~$x\in Q$, by inequality~\eqref{eqn:bloHaar}. So
\[
  \p_Q - \essinf_{x\in Q} \p(x)
  \le \sqrt{C_A} + C_B + C_3^d.
\]
Thus inequality~\eqref{eqn:A1} holds with constant~$C_3 =
\sqrt{C_A} + C_B + C_3^d$.

(b) In the case of $RH_\infty$, the same argument shows that
inequality~\eqref{eqn:RHinfinity} holds with constant~$C_4 =
C_4^d + C_B + \sqrt{C_A}$.
\end{proof}

\begin{proof}[Proof of Theorems~\ref{thm:Apexplog}
and~\ref{thm:RHpexplog}]
We prove Theorems~\ref{thm:Apexplog}
and~\ref{thm:RHpexplog} together.
Let $K$ denote any one of the classes $A_p$ with $1 \le p \le
\infty$, or $RH_p$ with $1 < p \le \infty$. For $\om$ in~$K$
let $K(\om)$ denote the corresponding $A_p$~constant or
$RH_p$~constant. Let $K^d$ denote the corresponding dyadic
class, and for $\om^t$ in $K^d$ let $K^d(\om^t)$ denote the
corresponding constant.

It follows from the definitions of~$A_p$ and~$RH_p$ that if
$\om$ is in~$K$ then for each constant~$\lambda > 0$ the
weight~$\lambda\om$ is also in~$K$, and $K(\lambda\om) =
K(\om)$.

For each $t\in[0,1]$ let $\om^t$ be a weight in~$K^d$, and let
$\p^t := \log\om^t$. Let $\p := \log\Om$. By hypothesis (iii)
of Theorem~\ref{thm:Apexplog} or~(iii$'$) of
Theorem~\ref{thm:RHpexplog}, the constants~$K^d(\om^t)$ are
bounded independently of~$t\in[0,1]$.

Without loss of generality, we may assume that for a.e.~$t\in [0,1]$ the $\p^t$
have mean value zero. To see this, note that the mean value $\p^t_\T :=
\intav_\T \p^t(x) \, dx$ is finite for a.e.~$t$ by hypothesis~(ii). Let
$\widetilde{\p}^t := \p^t - \p^t_\T$. Then for $x\in\T$,
\begin{align*}
    \Om(x)
    & = \exp\left\{\int_0^1 \p^t(x+t) \, dt\right\} \\
    & = \exp\left\{\int_0^1 \p^t_\T \, dt\right\}
        \exp\left\{\int_0^1 \widetilde{\p}^t(x+t) \, dt\right\} \\
    & = \exp\left\{\int_0^1 \!\!\! \int_0^1
        \log\om^t(x) \, dx \, dt\right\}
        \exp\left\{\int_0^1 \widetilde{\p}^t(x+t) \, dt\right\}.
\end{align*}
Denote the second term on the right-hand side
by~$\widetilde{\Om}(x)$. The first term on the right-hand side
is finite and nonzero by hypothesis~(ii) of Theorems~\ref{thm:Apexplog}
and~\ref{thm:RHpexplog}, and is positive. Thus $\Om$ is in $K$ if and only if
$\widetilde{\Om}$ is in~$K$. Moreover, $K(\Om) = K(\widetilde{\Om})$.

We show that hypotheses (i)--(iii) of
Lemma~\ref{lem:CACBestimates} hold here. The mapping
$t\mapsto\widetilde{\p}^t = (\log\om^t) - \p^t_\T$ is
measurable, since $t\mapsto\om^t$ is measurable by hypothesis.
By the dyadic version of Lemma~\ref{lem:AptoBMO}, the $\bmo_d$
norms $||\widetilde{\p}^t||_d = ||\p^t||_d$ are uniformly
bounded, by a constant depending only on~$p$ and on the bound
on the $K^d$~constants of the dyadic weights~$\om^t$. Each
$\widetilde{\p}^t$ has mean value zero, by construction.

For convenience we now drop the tildes, writing $\p^t$ for
$\widetilde{\p}^t$ and $\Om$ for $\widetilde{\Om}$ from here
on.

As an aside, we note that $\log\Om$ is in $\bmo$
(see~\cite{GJ}), and so by the John--Nirenberg Theorem the
function $\Om^\delta$ is in $A_p$ for $\delta > 0$ sufficiently
small. We now prove that $\Om$ itself is in~$A_p$.

\emph{Case~1:} $K = A_{\infty}$. By the dyadic version of
Lemma~\ref{lem:ApRHpcharn}(a), there is a constant $C_1^d$
depending only on the bound on the $A_\infty^d$~constants of
the~$\om^t$ such that each $\p^t$ satisfies
inequality~\eqref{eqn:betainequalitydyadic} with $\beta = 1$:
\[
    \intav_I \exp\{\p^t(x) - \p^t_I\} \, dx
    \le C_1^d.
\]
Take $C^d(1) = C_1^d$. Then by Lemma~\ref{lem:betainequality},
$\p = \log\Om$ satisfies inequality~\eqref{eqn:betainequality}
with $\beta = 1$: there is a constant $C(1)$ depending only
on~$C^d(1)$ such that
\[
    \intav_Q \exp\{\p(x) - \p_Q\} \, dx
    \le C(1).
\]
Take $C_1 = C(1)$. Lemma~\ref{lem:ApRHpcharn}(a) now implies
that $\Om\in A_{\infty}$, with $A_\infty$~constant bounded
by~$C_1$. The dependence of the constants is illustrated in the
upper row of Figure~\ref{fig:depconstsAp} (taking $p = \infty$
there). We see that $A_\infty(\Om)$ depends only on the bound
on the $A_\infty^d$ constants of the weights~$\om^t$.

\setlength{\unitlength}{1cm}

\begin{figure}[h]
    \begin{center}
        \begin{picture}(15,5)(0.5,0)
            \put(0,3){\framebox(2.4,1.4){$\underset{t\in[0,1]}{\max} \,
                                         A_p^d(\om^t)$}}
            \put(5.8,3){\framebox(2.7,1.4){\shortstack{$C_1^d = C^d(\beta)$,\\
                                                       with $\beta = 1$}}}
            \put(11.2,3){\framebox(1.2,1.4){$C_1$}}
            \put(14.4,1.5){\framebox(1.5,1.4){$A_p(\Om)$}}
            \put(0.6,0){\framebox(1.2,1.4){$p$}}
            \put(5.8,0){\framebox(2.7,1.4){\shortstack{$C_2^d = C^d(\beta)$,\\
                                               with $\beta = \frac{-1}{p-1}$}}}
            \put(11.2,0){\framebox(1.2,1.4){$C_2$}}
            \put(3,3.7){\vector(1,0){2.2}}
            \put(3,4){\shortstack{dyadic\\ Lemma~\ref{lem:ApRHpcharn}(a)}}
            \put(3,3){\vector(3,-2){2.2}}
            \put(3,0.7){\vector(1,0){2.2}}
            \put(3,1){\shortstack{dyadic\\ Lemma~\ref{lem:ApRHpcharn}(b)}}
            \put(9,3.7){\vector(1,0){1.8}}
            \put(9,4){\shortstack{Lemma~\ref{lem:betainequality},\\
                                  $\beta = 1$}}
            \put(9,0.7){\vector(1,0){1.8}}
            \put(9,1){\shortstack{Lemma~\ref{lem:betainequality},\\
                                  $\beta = \frac{-1}{p-1}$}}
            \put(12.8,3.7){\vector(3,-2){1.2}}
            \put(12.8,1.9){\shortstack{Lemma\\
                                  \ref{lem:ApRHpcharn}(b)}}
            \put(12.8,0.7){\vector(3,2){1.2}}
        \end{picture}
    \end{center}
    \caption{Dependence of the constants in the proof of Theorem~1,
        for the case $K = A_p$ with $1 < p < \infty$.}
    \label{fig:depconstsAp}
\end{figure}
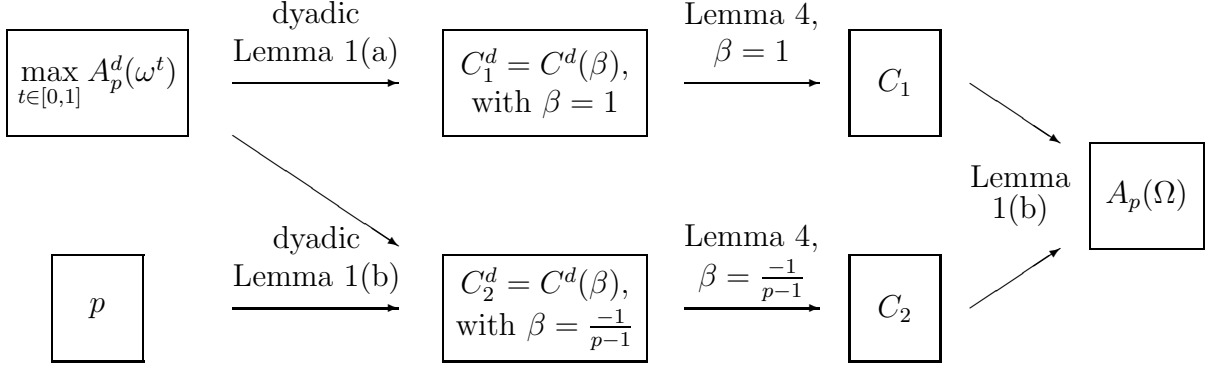

\emph{Case~2:} $K = A_p$, $1 < p < \infty$. As for case~1, using
Lemma~\ref{lem:ApRHpcharn}(b) and
Lemma~\ref{lem:betainequality} with both $\beta = 1$ and $\beta
= -1/(p-1)$. Figure~\ref{fig:depconstsAp} illustrates the
dependence of the constants. We find that $A_p(\Om)$ depends
only on~$p$ and on the bound on the $A_p^d$ constants of
the~$\om^t$.

\emph{Case~3:} $K = A_1$. As for case~1, using
Lemma~\ref{lem:ApRHpcharn}(c), Lemma~\ref{lem:betainequality}
with~$\beta = 1$, and Lemma~\ref{lem:A1RHinfinityinequality}.
The constant $A_1(\Om)$ depends only on the bound on the
$A_1^d$ constants of the~$\om^t$.

\emph{Case~4:} $K = RH_{\infty}$. As for case~1, using
Lemma~\ref{lem:ApRHpcharn}(d) and
Lemma~\ref{lem:A1RHinfinityinequality}. The constant
$RH_\infty(\Om)$ depends only on the bound on the $RH_\infty^d$
constants of the~$\om^t$.

\emph{Case~5:} $K = RH_p$. As for case~1, using
Lemma~\ref{lem:ApRHpcharn}(e) and
Lemma~\ref{lem:betainequality} with $\beta = p$. We find that
$RH_p(\Om)$ depends only on~$p$ and on the bound on the
$RH_p^d$ constants of the~$\om^t$.

This completes the proof of Theorems~\ref{thm:Apexplog}
and~\ref{thm:RHpexplog}.
\end{proof}

\emph{Remark.} An alternative proof of Theorem~\ref{thm:Apexplog} for $K =
A_p$, with $1 < p \leq \infty$, can be obtained as follows from the
result for $K = A_1$, using factorization of $A_p$
weights~\cite{Jon}.
Suppose $1 \le p < \infty$. If $\omega_1$ and $\omega_2$ are
$A_1$ weights, then $\omega = \omega_1 \omega_2^{1-p}$ is an
$A_p$ weight, with constant
$A_p(\om) \le A_1(\om_1) A_1(\om_2)^{p-1}$.
Conversely, if $\omega \in A_p$, then there exist $\omega_1$
and $\omega_2$ in $A_1$ such that $\omega = \omega_1
\omega_2^{1-p}$. The $A_1$ constants of~$\om_1$
and $\om_2$ depend only on~$p$ and on the $A_p$ constants
of~$\om$, as noted in~\cite[p.717]{Gra}. The analogous results
hold in the dyadic setting~$A_p^d$.


\begin{lem}\label{lem:pto1}
  If Theorem~\ref{thm:Apexplog} holds for $A_1$, then it holds for
  every~$A_p$, $1 \leq p \leq \infty$.
\end{lem}

\begin{proof}
By the (dyadic) factorization theorem, for each $t\in[0,1]$
there exist $\om^t_1$, $\om^t_2 \in A_1^d$ such that $\om^t =
\om^t_1(\om^t_2)^{1-p}$. Furthermore $A_1(\om_1^t)$ and
$A_1(\om_2^t)$ are uniformly bounded, independent of $t\in
[0,1]$, by a constant depending on~$p$ and on the bound for the
constants~$A_p(\om^t)$. Then
\begin{align*}
    \Om(x)
    &= \exp\left\{\int_0^1 \log\om^t(x+t) \, dt\right\} \\
    &= \exp\left\{\int_0^1
      \log\left[\om^t_1(x+t)\big(\om^t_1(x+t)\big)^{1-p}\right]
      \, dt\right\} \\
    &= \exp\left\{\int_0^1 \log\om^t_1(x+t) \, dt
      + (1-p) \int_0^1 \big(\om^t_1(x+t)\big) \, dt\right\} \\
    &= \left[\exp\left\{\int_0^1 \log\om_1^t(x+t) \, dt\right\}\right]
      \left[\exp\left\{\int_0^1 \log\om_2^t(x+t) \,
      dt\right\}\right]^{1-p}.
\end{align*}
By hypothesis, both of the expressions in square brackets are
in~$A_1$. Therefore $\Om$ is in~$A_p$ as required. Furthermore,
the $A_p$ constant of~$\Om$ depends only on the $A_p^d$
constants of the weights~$\om^t$.

The result for $A_\infty$ follows immediately from the result
for $A_p$ with $1 \leq p < \infty$, using the observations on
the dependence of the constants made before
Lemma~\ref{lem:ApRHpcharn} above.
\end{proof}

\section{Translation-averaging of doubling weights}
\label{sec:dbltranslproof}

It appears that the obstacle to the translation-average~$\om(x)
= \int_0^1 \om^t(x + t) \, dt$ of $A_p^d$ weights~$\om^t$ being
in~$A_p$ is that the assumption of the~$\om^t$ being dyadic
doubling is insufficient to guarantee that the
translation-average is actually doubling~\cite{W}. A natural
conjecture is that if the translation-average of a given family
of $A_p^d$~weights is in fact doubling, then this doubling
weight also belongs to true~$A_p$. As a step in this direction,
we show that the presumably stronger assumption that each of
the $A_p^d$ weights~$\om^t$ is doubling does imply that their
translation-average~$\om$ is in~$A_p$.

\begin{thm}\label{thm:doubling}
  Suppose $\{\om^t\}_{t\in [0,1]}$ is a family of doubling
  weights on~$\T$, with doubling constants bounded by a
  constant~$C_{{\rm dbl}}$ independent of~$t\in [0,1]$. Suppose
  the mapping $t \mapsto \om^t$ is measurable. Fix $p$ with $1
  \leq p \leq \infty$, and suppose each $\om^t$ is in $A_p^d$,
  with $A_p^d$~constant bounded by a constant~$V_{p,d}$
  independent of~$t\in [0,1]$. Then the translation-average
  \[
    \om(x) := \int_0^1 \om^t(x+t) \, dt
  \]
  belongs to~$A_p$ on~$\T$, with $A_p$ constant depending only
  on $p$, $V_{p,d}$, and~$C_{{\rm dbl}}$.

  Similarly, for $p$ with $1 < p \leq \infty$, if each~$\om^t$
  is in $RH_p^d$ and $RH_p^d(\om^t) \leq V_{p,d}$ for
  all~$t\in[0,1]$, then
  $\omega(x)$ belongs to~$RH_p$ on~$\T$, with $RH_p(\om)$
  depending on $p$, $V_{p,d}$, and $C_{{\rm dbl}}$.
\end{thm}

\begin{proof}
A doubling weight assigns comparable mass, with a constant
depending only on the doubling constant, to any given
interval~$Q$ and to each of the dyadic intervals at scale~$|Q|$
that intersect~$Q$. (This observation can fail if the weight is
dyadic doubling but not doubling.)

As a consequence, for fixed~$p$ with $1 \leq p \leq \infty$, an
$A_p^d$ weight is doubling if and only if it is actually
in~$A_p$. Moreover the $A_p$~constant depends only on the
$A_p^d$~constant and the doubling constant, which in turn
depend only on the $A_p$~constant. The same is true for~$RH_p$,
for $1 < p \leq \infty$. For example, one finds that for $A_p$
with $1 < p < \infty$, the $A_p$ constant of~$\om^t$ is bounded
by $V_p := 2^{2p - 1} \, V_{p,d} \, C_{\text{dbl}}^2$, while
for $RH_p$ with $1 < p < \infty$, the $RH_p$ constant
of~$\om^t$ is bounded by $V_p:= 2^{2/p} \, V_{p,d} \,
C_{\text{dbl}}^2$.

Theorem~\ref{thm:doubling} now follows easily for~$A_p$, $1 < p
< \infty$, using Muckenhoupt's original identification of~$A_p$
in terms of the boundedness of the Hardy--Littlewood maximal
function~$M$, defined as usual by $Mf(x) := \sup_{Q\ni x}
\intav_Q |f(y)| \, dy$. In particular, for these~$p$, a
nonnegative locally integrable function~$\om$ is in~$A_p$ if
and only if there is a constant~$C$ such that for all locally
integrable functions~$f$,
\begin{equation}\label{eqn:Apmaximalfn}
  \int_\R |Mf(x)|^p \, \om(x) \, dx
  \le C \int_\R |f(x)|^p \, \om(x) \, dx.
\end{equation}
Moreover, if inequality~\eqref{eqn:Apmaximalfn} holds then
$\om\in A_p$ and $A_p(\om) \le C$, while if $\om\in A_p$ then
inequality~\eqref{eqn:Apmaximalfn} holds with $C$ depending on
$p$ and~$A_p(\om)$.

%

For $RH_p$, $1 < p < \infty$, Theorem~\ref{thm:doubling}
follows from Minkowski's Integral Inequality and the
observation above on comparable mass. The cases of~$A_1$,
$A_\infty$, and~$RH_\infty$ are also straightforward, and we
omit the proofs.
\end{proof}


\section{Generalizations to the polydisc}\label{sec:polydisc}

We extend the above results for $A_p(\T)$ and $RH_p(\T)$ to the
setting of the polydisc. For ease of notation, the statements
and proofs given below are expressed for the bidisc. However,
they generalize immediately to the polydisc for arbitrarily
many factors.

The theory of product weights was developed by K.-C.~Lin in his
thesis~\cite{L}, while the dyadic theory was developed in
Buckley's paper~\cite{B}. The product $A_p$ and $RH_p$ weights
and the product doubling weights, and their dyadic analogues,
are defined exactly as in
Definitions~\ref{def:doubling}--\ref{def:RHp} in
Section~\ref{sec:definitions}, with intervals in~$\T$ being
replaced by rectangles in~$\T\otimes\T$. It follows that a
product weight belongs to $A_p(\T\otimes\T)$ if and only if it
belongs to $A_p(\T)$ in each variable separately.

To be precise, $\om\in A_p(\T\otimes\T)$ if and only if
$\om(\cdot,y) \in A_p(\T)$ uniformly for a.e.~$y\in\T$ and
$\om(x,\cdot) \in A_p(\T)$ uniformly for a.e.~$x\in\T$. In one
direction this is a consequence of the Lebesgue Differentiation
Theorem, letting one side of the rectangle shrink to a point.
The converse uses the equivalence between $\om\in
A_p(\T\otimes\T)$ and inequality~\eqref{eqn:Apmaximalfn} with
$M$ replaced by the strong maximal function~\cite[p.83]{S}.
Further, the $A_p(\T\otimes\T)$ constant depends only on the two
$A_p(\T)$ constants, and vice versa. The analogous
characterizations in terms of the separate variables hold for
product~$RH_p$ weights and for product doubling weights, and
for the dyadic product $A_p$, $RH_p$, and doubling weights.

\begin{thm}\label{thm:productApRHpexplog}
  Fix $p$ with $1 \le p \leq \infty$. Let
  $\{\om^{s,t}\}_{s,t\in [0,1]}$ be a family of dyadic~$A_p$
  weights on the boundary of the bidisc, $\om^{s,t}\in
  A_p^d(\T\otimes\T)$, such that
  \begin{enumerate}
  \item[{\rm (i)}] the mapping $(s,t)\mapsto\om^{s,t}$ is
      measurable,

  \item[{\rm (ii)}] the appropriate averages of the
      logarithms of the weights $\om^{s,t}$ are finite:
\[
      \int_0^1\!\!\!\int_0^1\!\!\! \int_0^1\!\!\!\int_0^1
            \big| \log \om^{s,t}(x,y)\big| \, dx \, dy \, ds \,dt
         < \infty,
\]
  and
  \item[{\rm (iii)}] the $A_p^d(\T\otimes\T)$~constants
      $A_p^d(\om^{s,t})$ are uniformly bounded, independent
      of $s$, $t\in [0,1]$.
  \end{enumerate}
  Then the product geometric-arithmetic average
  \[
    \Om(x,y)
    := \exp\left\{\int_0^1 \!\!\! \int_0^1
	\log \om^{s,t}(x+s,y+t) \, ds \, dt\right\}
  \]
  lies in $A_p(\T\otimes\T)$, with $A_p$ constant depending on~$p$ and
  on the bound on the $A_p^d$ constants of the
  dyadic weights~$\om^{s,t}$.

  Similarly, for fixed $p$ with $1 < p \leq \infty$, if each
  $\om^{s,t}$ is in $RH_p^d(\T\otimes\T)$, if hypotheses~\textup{(i)}
  and~\textup{(ii)} hold, and if
  \begin{enumerate}
    \item[{\rm (iii$'$)}] the
        $RH_p^d(\T\otimes\T)$~constants $RH_p^d(\om^{s,t})$
        are uniformly bounded, independent of $s$, $t\in
        [0,1]$,
  \end{enumerate}
  then $\Om(x,y)$ lies in $RH_p(\T\otimes\T)$, with $RH_p$~constant
  depending on~$p$ and on the bound on the $RH_p^d$ constants of
  the~$\om^{s,t}$.
\end{thm}

\begin{proof}
The proof is by iteration of the one-parameter argument,
relying on Lemma~\ref{lem:ApRHpcharn} and Theorems~\ref{thm:Apexplog}
and~\ref{thm:RHpexplog}. We give the argument for $K = A_p$, $1 < p < \infty$.
The other cases follow similar iteration arguments and we omit the proofs.

We show that, for a.e.\ fixed $y$, $\Omega(x,y)$ belongs to $A_p(\T)$ in
the variable~$x$, with constant independent of~$y$.
The hypotheses of Theorem~\ref{thm:Apexplog}
follow immediately from our assumptions. In particular,
$s \mapsto \omega^{s,t}(\cdot,y)$ is measurable for each $t$ and $y$,
and $\om^{s,t}(\cdot,y)$ belongs to $A_p^d(\T)$ in the first variable for
all $s$, $t$ and for a.e.\ $y$, with constants independent of $s$, $t$, and $y$.
Fix such a $y$; for emphasis we'll denote it by~$y^*$.

The conclusion of Theorem~\ref{thm:Apexplog} is that the function
\[
\Omega_1(x,y^*)
:= \exp\left\{\int_0^1 \log\om^{s,t}(x+s, y^*)\,ds\right\}
\]
belongs to $A_p(\T)$ in $x$, with constant independent of~$y^*$. Let
\[
    \p_1(x,y^*)
    := \log \Om_1(x,y^*)
    = \int_0^1 \log \om^{s,t}(x+s,y^*) \, ds.
\]
By Lemma~\ref{lem:ApRHpcharn}(b) there are constants~$C_1$ and $C_2$, depending
only on $p$ and on the $A_p^d(\T\otimes\T)$ constants of the~$\om^{s,t}$, such
that
\begin{align}\label{eqn:phi1Apineq1}
    \sup_Q \intav_Q \exp\left\{\p_1(x,y^*)
        -  \big(\p_1(\cdot,y^*)\big)_Q\right\} \, dx
    &\leq C_1, \\
    \label{eqn:phi1Apineq2}
    \sup_Q \intav_Q \exp\left\{\frac{-1}{p-1}
        \Big[\p_1(x,y^*) - \big(\p_1(\cdot,y^*)\big)_Q\Big]\right\} \, dx
    &\leq C_2.
\end{align}
We show that the same inequalities hold when $\p_1$ is replaced by
\[
    \p(x,y^*)
    := \log \Om(x,y^*)
    = \int_0^1 \!\!\!\int_0^1 \log \om^{s,t}(x+s,y^*+t) \, ds \, dt.
\]

Fix an interval~$Q$. An application of Fubini's Theorem shows that
\begin{align*}
    \lefteqn{\p(x,y^*) - \big(\p(\cdot,y^*)\big)_Q} \hspace{1cm} \\
    &= \int_0^1\!\!\!\int_0^1 \p^{s,t}(x+s,y^*+t) \, ds \, dt
        - \intav_Q\!\int_0^1\!\!\!\int_0^1 \p^{s,t}(x'+s,y^*+t)
        \, ds \, dt \, dx' \\
    & = \int_0^1 \Big[\p_1(x,y^*+t) - \big(\p_1(\cdot,y^*+t)\big)_Q\Big] \, dt.
\end{align*}
Then by Jensen's inequality and Tonelli's Theorem,
\begin{align*}
    \lefteqn{\intav_Q \exp\left\{\p(x,y^*)
        - \big(\p(\cdot,y^*)\big)_Q\right\} \, dx} \hspace{1cm} \\
    &\leq \intav_Q \!\int_0^1 \exp\left\{\p_1(x,y^*+t)
        - \big(\p_1(\cdot,y^*+t)\big)_Q\right\} \, dt \, dx \\
    &= \int_0^1\!\!\intav_Q \exp\left\{\p_1(x,y^*+t)
        - \big(\p_1(\cdot,y^*+t)\big)_Q\right\} \, dx \, dt \\
    &\leq C_1,
\end{align*}
by inequality~\eqref{eqn:phi1Apineq1} with $y^*$ replaced by $y^* + t$ for
a.e.~$t$.

The same argument can be used to verify that inequality~\eqref{eqn:phi1Apineq2}
holds for~$\p(\cdot,y^*)$ for a.e.~$y^*$. Therefore by
Lemma~\ref{lem:ApRHpcharn}, $\Om(x,y^*)$ belongs to $A_p(\T)$ in~$x$ for
a.e.~$y^*$, with uniform constants.

In an identical fashion, we find that $\Om(x^*,y)$ belongs to $A_p(\T)$ in~$y$
for a.e.~$x^*$, with uniform constants, which proves the theorem for $K = A_p$,
$1 < p < \infty$.
\end{proof}

\emph{Remark.} As in the one-parameter case, there is an alternative proof of
the geometric-arithmetic averaging result (Theorem~\ref{thm:productApRHpexplog})
for $A_p(\T\otimes\T)$ where $1 < p \le \infty$, relying on the
$A_1(\T\otimes\T)$ case and the generalization to the bidisc setting~\cite{Jaw}
of the $A_p$ factorization theorem. Moreover, the product $A_p$ case can also be
derived from the one-parameter result using the maximal function
characterization of this weight class.

For product doubling weights, we also have a result analogous to
Theorem~\ref{thm:doubling}.

\begin{thm}\label{thm:Bidoubling}
  Suppose $\{\om^{s,t}\}_{s,t\in [0,1]}$ is a family of
  doubling weights on~$\T\otimes\T$, with doubling constants
  bounded by a constant~$C_{{\rm dbl}}$ independent of~$s$,
  $t\in [0,1]$. Suppose the mapping $(s,t) \mapsto \om^{s,t}$ is
  measurable. Fix $p$ with $1 \leq p \leq \infty$, and suppose
  each $\om^{s,t}$ is in $A_p^d(\T\otimes\T)$ with
  $A_p^d(\T\otimes\T)$~constant bounded by a constant~$V_{p,d}$
  independent of~$s,t\in [0,1]$. Then the translation-average
  \[
    \om(x,y)
    := \int_0^1 \!\!\! \int_0^1
      \om^{s,t}(x+s,y+t) \, ds \, dt
  \]
  belongs to~$A_p(\T\otimes\T)$, with an $A_p$ constant
  depending only on $p$, $V_{p,d}$, and $C_{{\rm dbl}}$.

  Similarly, for $1 < p \leq \infty$, if the $A_p^d(\T\otimes\T)$ assumption
  above is replaced by the assumption that each $\om^{s,t}$ is
  in~$RH_p(\T\otimes\T)$ with $RH_p^d(\om^{s,t}) \leq V_{p,d}$ for all~$s$,
  $t\in [0,1]$, then $\omega(x,y)$ belongs to~$RH_p(\T\otimes\T)$, with
  $RH_p(\om)$ depending on $p$, $V_{p,d}$, and~$C_{{\rm dbl}}$.
\end{thm}

The proof is by iteration of the one-parameter argument given
for Theorem~\ref{thm:doubling} above.

\end{document}